# Alternative ways to initiate students' intuition, and hence internalization, of why zero factorial is equal to one

**Munir Mahmood[1], Lori Murray[2] , Ricardas Zitikis[3], Ibtihal Mahmood[4]**

[1]Mathematics and Statistics Education Australia, Unit 1, 39 Moodemere Street, Noble Park, Victoria 3174, Australia
[2]School of Management, Economics, and Mathematics, King's University College at Western University, London, Ontario, N6A 2M3, Canada
[3]School of Mathematical and Statistical Sciences, Western University, London, Ontario, N6A 5B7, Canada
[4]University of Melbourne, Parkville 3052, Melbourne, Australia

*Correspondence: Lori Murray, lori.murray@uwo.ca

## Abstract

Zero factorial, defined to be one, is often counterintuitive to students but nonetheless an interesting concept to convey in a classroom environment. The challenge is to delineate the concept in a simple and effective way through the practice of justification, a familiar concept in mathematics and science education. In this regard two algebraic and one statistical justification using the squeeze theorem are presented. To assess the effectiveness of the justifications, a student survey was conducted at a comprehensive university incorporating the analysis of the pre- and post-presentation statements. They clearly present that the justifications are useful in giving credence to zero factorial equals one. Overall, the results from the online survey supports that the students preferred Justification 1. The justifications provide instructors alternative ways to initiate exploration of students' intuitive set up of comprehending unobvious facts like zero factorial equals one.

For a range of learners with their varied abilities to perform various mental activities most closely associated with learning and problem-solving, the justifications as simple alternative methods offer the potential to raise the current level of cognitive skills to inspire differentiated paths of learning. These are evident from survey results noting the role of statistical thinking and techniques.

**Keywords**: bounds, factorial, gamma function, justification, squeeze theorem

## INTRODUCTION

Improving numeracy at all levels of education has been a challenging task (Aunio et al., 2021; Brumwell & MacFarlane, 2020; Galligan, 2013). Naturally, different pedagogical methods and techniques are needed and therefore have been employed at different levels (primary, high, university, etc.) and in different cultures, to make the explanations relevant to students. There is a consensus (see, e.g. (Risdiyanti et al., 2024)) that Realistic Mathematics Education (RME) approach and its Indonesian adaptation, known as Pendidikan Matematika Realistik Indonesia (PMRI), have emerged as promising ways to facilitate student learning. As an illustrative specimen of how this can be achieved, we note such innovations as the use of traditional food (Gustiningsi et al., 2023), dance (Rawani et al., 2023), and nature tourism (Putra et al., 2023) for teaching various mathematical concepts. For many other enlightening illustrations, we refer to the proceedings of NACOME (2023).

To illustrate how challenging it is to teach mathematical subjects, let us look at the factorial, which plays a prominent role in areas such as mathematics, statistics, and probability. To see the challenge, we first recall that by its very definition, the factorial of any natural number is defined as the product of that number until we reach 1. Hence with this descending trend of factorial calculation explanation, in the case of calculating 0 factorial, the same definition would intuitively result in $0$. This invariably leads to contradict many proven mathematical results. So, to alleviate this issue, zero factorial is then defined to be one. This definition, however, is often counterintuitive to students.

Often in high school it is taught that zero factorial is one by convention, which only adds further void to learning of this already unobvious fact. At the university level, quite often in teaching introductory courses, students are merely instructed that zero factorial equals one as a fact for convenience and are not provided a methodical explanation as to why (Buglear & Castell, 2019). At the advanced university level, factorials are generalized into the gamma function, thus implying that zero factorial is one by a fancy integral associated with the gamma function, which is apparently an attempt to mend the conceptual learning from high school. Unfortunately, many students seem to accept it unknowingly without asking for any in depth or inherent justification of it because by then they have already learned integration from a first-year calculus course. The integration works out well and so it is deemed to be true with no further question asked. In this article we concentrate on teaching why zero factorial equals one in introductory mathematics, statistics, and other related courses.

Understanding, intuitively appreciating, and internalizing why zero factorial equals one is important for several reasons. It is a fundamental concept in combinatorics and algebra and serves as a basis for many mathematical principles and calculations. In combinatorics, factorial plays a crucial role in counting permutations and combinations. Understanding the value of zero factorial ensures that combinatorial formulas are correctly applied and interpreted. Many algorithms and programs involve factorial calculations, and understanding the value of zero factorial is essential for writing correct and efficient codes in computer science applications. Exploring why zero factorial equals one deepens on one's understanding of mathematical concepts such as recursion, limits, and mathematical induction. Overall, understanding why zero factorial equals one is essential for a thorough comprehension of mathematics and its applications across various fields.

The connection of zero factorial to the definition of the gamma function, which with the help of computers can be visualized, provides a first-hand conceptual understanding of the concept of zero factorial. But for the purpose of teaching, the gamma function is not particularly helpful from the pedagogical point of view in early years of study, as it is quite challenging to explain the rationale behind the origin of the integral that defines the gamma function.

There are many ad hoc explanations of why zero factorial equals one, and arguments ensue whether, and to what extent, they can be considered mathematically accurate. Obviously, for the sake of facilitating intuition, at least during the first stages of explanation, we are bound to slightly depart from rigor, which could later be, and is, re-stored via the celebrated gamma function in more advanced calculus courses. It is a path that, in some way, resembles the use of Venn diagrams which hardly ever give precise descriptions of events, although provide most intuitively appreciated ways to explain what has been proved, or will be proved.

Coming back to the factorial, it should be noted at this point, however, that there is no simple way to express the factorial function $h(n) = n!$ beyond the positive integers $n$, and thus to reveal its value at $n = 0$. In particular, for the sole purpose of conveying a meaningful concept of the zero factorial, a simpler approximation to the integral form of the gamma function to date is unavailable. Hence it is desirable to derive or present a formula that would transition the factorial function for positive integer values to non-integer values, particularly to those that are in the interval (0,1), and then approach 0 along the curve.

To mitigate the above learning constraints associated with the gamma function, we present three novel intuitive pedagogical justifications, all of which providing visualizations and thus facilitating internalization. They all rely on the squeeze theorem to convince that zero factorial (i.e., 0!) is equal to 1, and not to 0 as many students would expect. The core concept of these justifications draws its strength from defining a symbolic $t!$ (read as $t$ factorial) in $t \in (0, 1)$ where the squeeze theorem remains valid but by interchanging their bounds in comparison to their bounds when $n$ is a positive integer greater than or equal to 1.

It may be noted that justification, as a practice, is at the heart of mathematics and may be regarded as one of the heuristic best practices, particularly when student learning and enhancing of their understanding of mathematics concepts are involved. In this regard Staples et al (2012) give an excellent account of justification as a practice in providing several ways to reach every student, as it offers and enables differentiated paths of learning. Having several choices is also useful for the instructor in selecting a most suitable one for the level and perceived inclinations of a given group of students (Dare et al., 2021). We illustrate how data and their statistical analyses can facilitate the selection of an effective justification based on a survey conducted at a comprehensive university.

In general, surveys have played a particularly important role in educational decision making, and one of the most researched databases is the OECD Education Statistics. On a smaller, localized scale, educational researchers have conducted surveys on issues such as the impact of active learning in undergraduate science, technology, engineering, mathematics, and medicine (STEMM) education (Wiggins et al., 2017), assessment design (Murray & Wilson, 2021), students' views on introductory science, technology, engineering, and mathematics (STEM) courses (Meaders et al., 2020), and experiential learning perspectives (Roberts et al., 2018; Zhai et al., 2017).

A well-designed survey is crucial in achieving statistically sound and trustworthy results (Bellhouse, 1984; Thompson, 1997; Wu & Thompson, 2020). The focus of the survey in the present paper is quite distinct: we aim not at general instructional matters but at a very specific mathematical concept and, given several alternatives, at choosing a most accessible way to present the concept in a classroom environment.

Apart from the foregoing aim of survey, quite often in introductory courses, instructors emphasize learning concepts by means of repetition and ample practice problems. This is arguably carried out to reinforce the initiation of students' intuition and facilitate internalization of concepts to make them easier to remember by building a learning process by associating what needs to be learned with what has already been learned. For this task, the development of curricula plays a decisive role. In this regard, Chen and Zitikis (2017, 2020) have explored the effects of learning a subject due to prior subjects learned, based on real data and in-depth statistical analyses.

We believe, in essence, that repeated teaching of at least one of the three justifications proposed, i.e., the one which stands out to be effective from repeated survey results, would enable nurture students' intuition to internalize and assimilate conceptually that zero factorial is indeed one, and not to be accepted only as a rigorous definition.

As we have noted earlier, in more advanced courses, this goal is effectively achieved via the well-known gamma function, which is $\Gamma(t) = \int_0^\infty x^{t-1}e^{-x}\,dx$ for all real $t > 0$. For details, see Devore (2012). The relationship $\Gamma(n+1) = n\Gamma(n) = n!$ between the gamma and factorial functions, and the equation $\Gamma(t+1) = t\Gamma(t)$ that holds for all $t \in (0,1)$, create a genuine bridge between the cases $n = 1$ and $n = 0$. Intuitively, if we take the limit as $t \downarrow 0$ inside the integral $\int_0^\infty x^t e^{-x}\,dx$, and all limits in this paper approach 0 from the right, we are left with the exponential function, which integrates to 1. This shows that $\Gamma(t+1)$ converges to 1 when $t \downarrow 0$. Defining 0! as the limit of $\Gamma(t+1)$ when $t \downarrow 0$, we can now claim that 0! is indeed equal to 1. This explanation supports the presentations of $0! = 1$ in Mahmood & Mahmood (2015) and Mahmood & Romdhane (2015), although it is not explicitly mentioned in those articles.

We aim to seek further simplified algebraic and statistical justifications and then let survey data choose the "best" one among them. The "best" is of course contingent on the university and the type of students. Given the prerequisites for enrolling into courses, the notion of "best" would then carry on for a given course for several years after the survey results. Such availability of the choicest justification addresses and serves the essential purpose of effective understanding of zero factorial in classes that may not have yet encountered the concept of the gamma function, or may have seen it in other contexts such as statistics but have not associated it with real numbers other than positive integers.

In this regard, we acknowledge to have greatly benefited from historical notes of Hairer and Wanner (1996) on the gamma function, factorials, and Euler's interest and achievements in interpolating factorials to non-integer values. Of course, the use of the gamma function would give the most mathematically rigorous justification of 0!=1, whereas Mahmood & Mahmood (2015) and Mahmood & Romdhane (2015) presented a particularly intuitive explanation.

The justifications in the current paper slightly reduce the intuitive appeal but add more mathematical rigor to the justification of 0!=1, hence exhibiting the delicate balance that we desire to achieve in a classroom setting when explaining new concepts to students. The issue is that the presence of the integral in the definition of the gamma function presents a tricky form, one which is not easy to absorb by students at such an early exposure stage of the concepts of $n$ factorial and zero factorial. Moreover, it is hard to explain the rationale behind the basis of such an integral in $t$ and $x$, of which $t$ ascertains the factorial calculation while the embedded exponential function assures that the definite integral produces any given factorial value.

To facilitate our following arguments in a gentle and intuitive fashion, a symbolic $t!$ is introduced in the next section, which allows a movement from the discrete $n$ to the continuous $t \in (0, 1)$. This enables taking the limits in a subtle way, aligning more with the gamma function but avoiding its extreme mathematical precision and thus complexity. The justifications provide both a mathematical and visual explanations as to why zero factorial equals one, and they are either in addition to Mahmood & Mahmood (2015) and Mahmood & Romdhane (2015) or improved versions of them.

# METHODS

## Mathematical Justifications

As explained at the end of the previous section, a symbolic $t!$ is applied to the three justifications with $t \in (0, 1)$. This provides the mechanism to distinguish symbolic $t!$ from the definition of the factorial $n$, which is valid when $n$ is a positive integer. The factorial of $n$ is defined by

$$n! = n \cdot (n-1) \cdot \ldots \cdot 3 \cdot 2 \cdot 1 \qquad (1)$$

for integer $n \geq 1$, but for $n = 0$. Zero factorial, symbolically written as $0!$, is defined as $0! = 1$. In this section we present some easier justifications of $0! = 1$, and which of them is the "best" will then be decided by a data-driven statistical analysis in the results section. As already noted earlier, the three justifications will rely on the squeeze theorem which is restated below.

The Squeeze Theorem (Adams & Essex, 2010): If for three functions $f$, $g$, and $h$ we have the relationship $f(t) \leq g(t) \leq h(t)$ for all $t$ in some open interval containing $t_0$, except possibly at $t = t_0$, and if $f(t)$ and $h(t)$ have the same limit, say $L$, when $t$ approaches $t_0$, then $g(t)$ also has the same limit $L$ when $t$ approaches $t_0$.

In the following justifications, we shall have various functions $f$, $g$, and $h$, with $t_0 = 0$ and $L = 1$ in all of them.

For Justification 1, restating from Mahmood & Mahmood (2015), a lower bound for $n!$ is given by
$1 = 1^n = 1 \cdot 1 \cdot \ldots \cdot 1 \cdot 1 \cdot 1 \leq n!$ whereas we now propose using the following lower bound

$$\left(\frac{n}{2}\right)^{n/2} \leq n!. \qquad (2)$$

We can prove bound (2) as follows. Let $k$ be the (unique) positive integer such that $k - 1 \leq \frac{n}{2} < k$. Then

$$
\begin{aligned}
n! &= n \cdot \ldots \cdot k \cdot \ldots \cdot 2 \cdot 1 \\
&\geq n \cdot \ldots \cdot k \\
&\geq k^{n-k+1} \\
&\geq \left(\frac{n}{2}\right)^{n/2}
\end{aligned}
$$

where the right-most inequality holds because $n - k + 1 \geq \frac{n}{2}$ due to $\frac{n}{2} \geq k - 1$, with $k > \frac{n}{2}$. This completes the proof of bound (2).

Further, it was noted that an upper bound for $n!$ follows from

$$n! \leq n \cdot n \cdot \ldots \cdot n \cdot n \cdot n = n^n. \tag{3}$$

To avoid the application of l'Hospital's Rule to demonstrate $0! = 1$ of Mahmood & Mahmood (2015), we next propose the following upper bound

$$n^n \leq 2^{n^2}, \tag{4}$$

which is valid since $n < 2^n$. The information conveyed in statements (1) to (4) establishes the inequalities

$$\left(\frac{n}{2}\right)^{n/2} \leq n! \leq 2^{n^2}. \tag{5}$$

These inequalities are correct for positive integers, but we wish to see where they take us when $n$ goes to 0. Recall, however, that $n$ is a positive integer, and thus limits such as $n \to 0$ do not make sense. Hence, we build a bridge via a parameter $t \in (0, 1)$ and a symbolic $t!$, which is merely a function of $t$ squeezed between $\left(\frac{t}{2}\right)^{t/2}$ and $2^{t^2}$. Thus, we write the inequalities

$$\left(\frac{t}{2}\right)^{t/2} \leq t! \leq 2^{t^2}. \tag{6}$$

The left- and right-hand sides of (6) converge to 1 when $t \downarrow 0$. Hence as a consequence of the squeeze (or "two policemen and a drunk") theorem, defining $0!$ as the limit of $t!$ when $t \downarrow 0$, we conclude that $0! = 1$. Figure 1 visualizes these arguments by depicting the lower bound $a(t) = \left(\frac{t}{2}\right)^{t/2}$, the upper bound $b(t) = 2^{t^2}$, and the symbolic $t! \equiv \mathrm{factorial}(t)$ drawn in a wiggly way that we would use to depict on the whiteboard in a classroom set up when $t \in (0, 1)$.

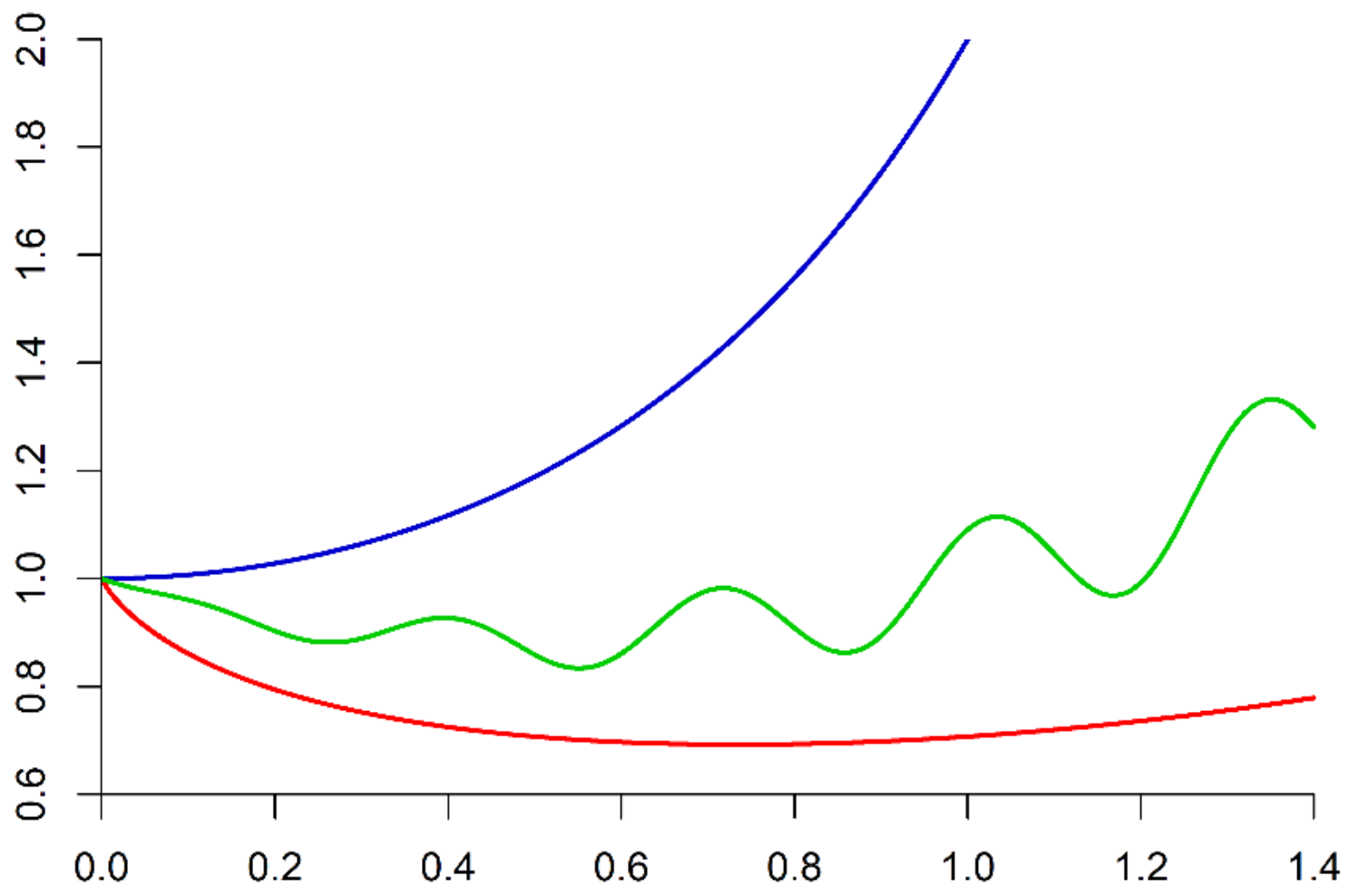

**Figure 1.** Pictorial summary of Justification 1 with the curves $a(t)$ (red), $b(t)$ (blue), and symbolic $t!$ (green).

For Justification 2, we bypass the rigorous presentation of using l'Hospital's Rule and the harmonic mean (*HM*) of Devore (2012) in order to arrive at zero factorial to be one. Applying the well-known inequality $GM \leq AM$ to the set $\{1,2,3,\ldots,n\}$, where *GM* and *AM* are the geometric and arithmetic means, gives the bound

$$(1 \cdot 2 \cdot \ldots \cdot n)^{1/n} \leq \frac{1+2+\ldots+n}{n}. \tag{7}$$

Expressing the left-hand term of (7) by the $n$ factorial notation given by (1) and also using the equation $1 + 2 + \ldots + n = \frac{n(n+1)}{2}$, it is clear that

$$1 \leq n! \leq \left(\frac{n+1}{2}\right)^n. \tag{8}$$

Again, these bounds are valid for positive integers, but we wish to have a bridge that connects $n = 1$ with $n = 0$, and thus we again rely on the symbolic $t!$ for all $t \in (0,1)$. Accordingly, we write the inequalities

$$\left(\frac{t+1}{2}\right)^t \leq t! \leq 1. \tag{9}$$

Remark 1: The inequalities in (8) are true for $n \geq 1$, but when $t$ is in (0, 1), the two "policemen" interchange their sides as noted by the inequalities in (9), keeping the "drunk", which is $t!$, still between themselves.

The left-hand side of (9) converges to 1 when $t \downarrow 0$, and so by the squeeze theorem and following the definition of $0!$ as the limit of $t!$ when $t \downarrow 0$, we conclude the equation $0! = 1$ from (9).

Figure 2 visualizes the above arguments by presenting $a(t) = \left(\frac{t+1}{2}\right)^t$, $b(t) = 1$, and the symbolic $t!$ for $t \in (0,1)$.

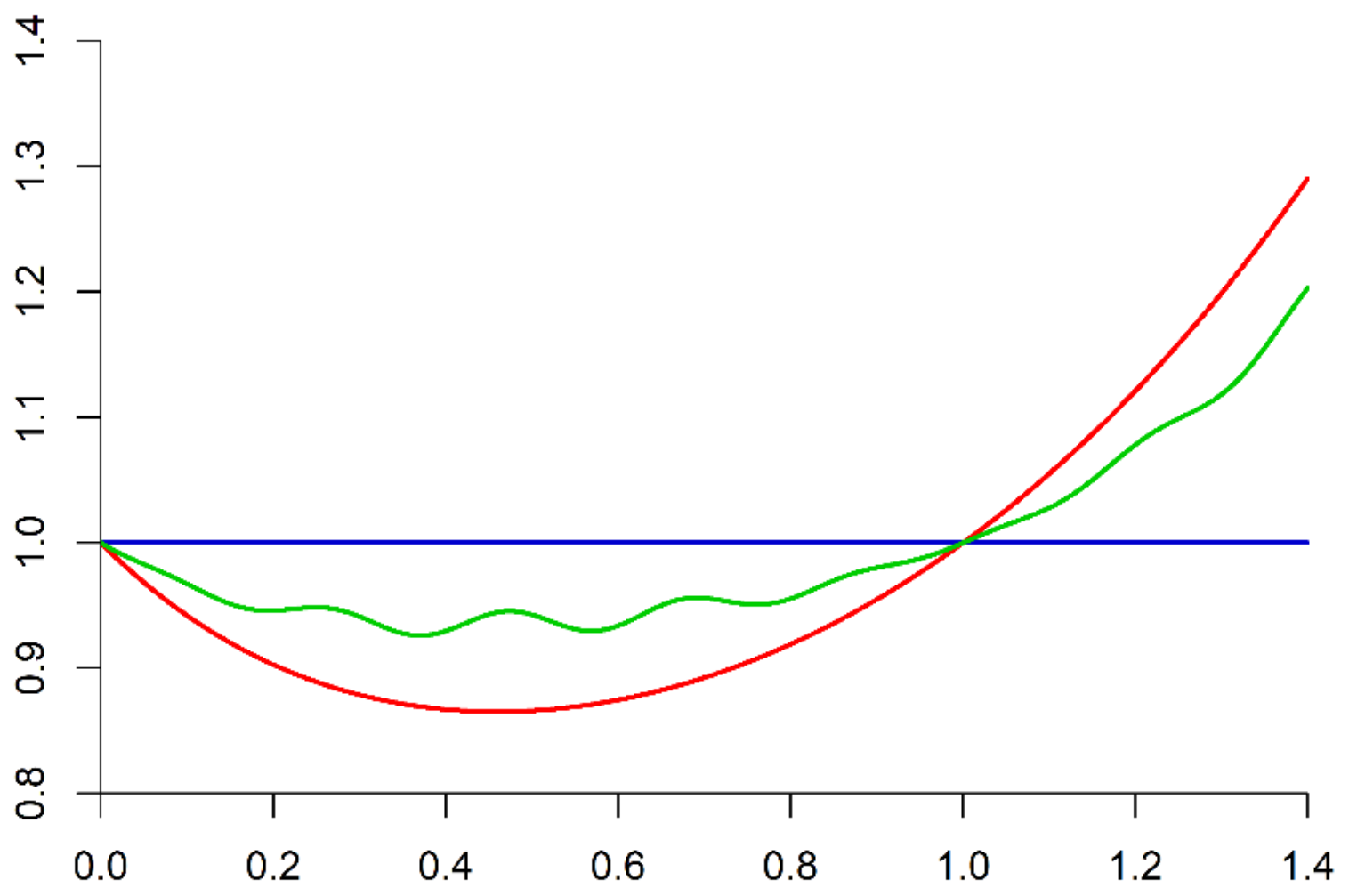

**Figure 2.** Pictorial summary of Justification 2 with the curves $a(t)$ (red), $b(t)$ (blue), and symbolic $t!$ (green).

For Justification 3, we consider a random variable $X$ that follows the exponential distribution with rate $\lambda = 1$, that is, its mean is $E[X] = 1$ and the probability density function is given by

$$f(x) = e^{-x} \tag{10}$$

for $x > 0$. For details of the exponential distribution, see (Devore, 2012). The *n*-th moment of $X$ is

$$E[X^n] = \int_0^\infty x^n f(x)\, dx = \int_0^\infty x^n e^{-x}\, dx = n!, \tag{11}$$

that is, $n! = E[X^n]$.

We now again create a bridge between $n = 1$ and $n = 0$ by introducing bounds $a(t)$ and $b(t)$ such that $a(t) \leq E[X^t] \leq b(t)$ for $t \in (0, 1)$. We start with an upper bound. For any concave function $g(x)$, the Jensen's inequality states that

$$E[g(X)] \leq g[E(X)]. \tag{12}$$

Applying this inequality for the concave function $g(x) = x^t$ for any $t \in (0, 1)$ yields

$$E[X^t] \leq [E(X)]^t = 1^t = 1. \tag{13}$$

Therefore, $b(t) = 1$. We next find a lower bound $a(t)$. We have

$$E[X^t] = \int_0^\infty e^{-x^{1/t}} dx \geq \int_0^1 e^{-x^{1/t}} dx \geq \int_0^1 \left(1 - x^{1/t}\right) dx$$

because $e^{-\theta} \geq 1 - \theta$ for every $\theta \geq 0$, and in our case $\theta = x^{1/t}$. Hence,

$$E[X^t] \geq 1 - \int_0^1 x^{1/t} dx = 1 - \frac{1}{1+\frac{1}{t}} = \frac{1}{1+t}$$

and so $a(t) = \frac{1}{1+t}$. In summary, we have

$$\frac{1}{1+t} \leq E[X^t] \leq 1 \tag{14}$$

and thus whatever the graph of symbolic $t!$ we may wish to draw on the classroom whiteboard, it must satisfy the bounds

$$\frac{1}{1+t} \leq t! \leq 1. \quad (15)$$

Since the left-hand side of bounds (15) converges to 1 when $t \downarrow 0$, we conclude that the limit of $t!$ when $t \downarrow 0$, which is by definition equal to $0!$, must be equal to 1.

Figure 3 visualizes $a(t) = \frac{1}{1+t}$, $b(t) = 1$, and the symbolic $t!$ in a wiggly way that we would use to draw on the classroom whiteboard for $t \in (0, 1)$.

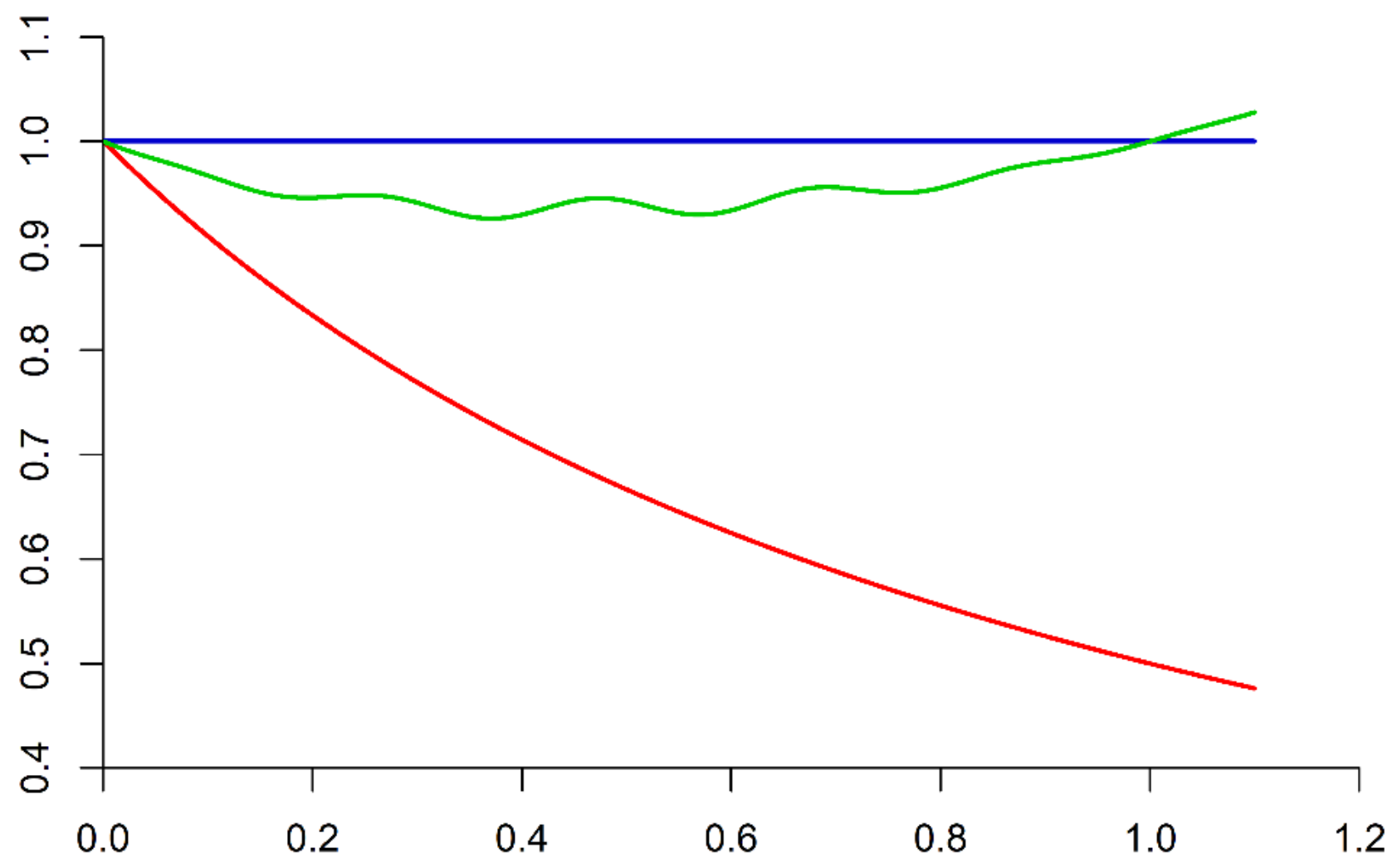

**Figure 3.** Pictorial summary of Justification 3 with the curves $a(t)$ (red), $b(t)$ (blue), and symbolic $t!$ (green).

## Survey Design and Participants

To assess the effectiveness of the three justifications pedagogically, a class presentation and an online survey were implemented. Our study sample comprised sixty-four undergraduate students enrolled in an introductory statistics course at a Canadian university. These students were pursuing majors in mathematics, finance, and economics and were invited to participate in the study. The student survey encompassed an 11-statement questionnaire utilizing a 5-point Likert-type scale ranging from 1 (strongly disagree) to 5 (strongly agree) administered via Qualtrics software. Approval for this study was obtained from the university institutional ethics review board.

The course was conducted online, utilizing synchronous online lectures for delivery. Consequently, the presentation took place online following one of these lectures. Prior to the presentation, the instructor prompted students to provide a rating for the initial statement (Statement 0) regarding the believability of zero factorial equals one. Once students submitted their ratings for Statement 0 via the Qualtrics software, they were presented with the three justifications, both verbally and visually as described above. Following the presentation, students were encouraged to provide ratings for the remaining 10 statements (Statements 1a to 4) via the online survey. The 11-statement questionnaire is outlined below:

0. *Zero factorial equals one (i.e., $0! = 1$) is believable.*

1a. *The explanation of symbolic t! in Justification 1 is helpful.*

| | |
|---|---|
| 1b. | *The pictorial summary of Justification 1 by the graph on slide 5 aided further to conceptualize 0! to be 1.* |
| 1c. | *Justification 1 is helpful in believing that* $0! = 1$ |
| 2a. | *The explanation of symbolic t! in Justification 2 is helpful.* |
| 2b. | *The pictorial summary of Justification 2 by the graph on slide 5 aided further to conceptualize 0! to be 1.* |
| 2c. | *Justification 2 is helpful in believing that* $0! = 1$ |
| 3a. | *The explanation of symbolic t! in Justification 3 is helpful.* |
| 3b. | *The pictorial summary of Justification 3 by the graph on slide 5 aided further to conceptualize 0! to be 1.* |
| 3c. | *Justification 3 is helpful in believing that* $0! = 1$ |
| 4. | *Overall, the statement that zero factorial equals one (i.e.,* $0! = 1$*) is more believable after having read the justifications.* |

### Data Analysis

The students' individual responses for each statement were assigned numerical values: strongly agree =5, somewhat agree =4, neither agree nor disagree =3, somewhat disagree =2, and strongly disagree =1. To assess the students' perceptions and to see which of the justifications were preferred, the aggregate counts and percentages of the ratings for each of the 11 statements were analyzed. In addition, a matched-pairs t-test was conducted to see if there was an improvement in the responses between the pre- and post-presentation statements.

## RESULTS AND DISCUSSION

### Survey Results

The survey data collected indicated that 62 students responded to all of the statements including the pre- and post-statements of the questionnaire. A summary of the data gathered from the online survey are shown in Tables 1 to 4.

Table 1 shows the counts and percentages for the statement "the explanation of symbolic t! is helpful" for each justification. The aggregate counts for the categories "strongly agree" and "somewhat agree" are 55/62 (88.71%), 51/62 (82.26%), and 45/62 (72.85%) for the three justifications, respectively, and indicate that the explanation of symbolic t! is the most helpful in Justification 1.

Table 2 shows the counts and percentages for the statement "the pictorial summary by the graph aided further to conceptualize 0! to be 1" for each justification. The aggregate counts for "strongly agree" and "somewhat agree" are 57/62 (91.94%), 57/62 (91.94%), and 51/62 (82.26%) for the three justifications, respectively, and indicate that the graphs used in Justifications 1 and 2 are the most helpful.

Table 3 shows the counts and percentages for the statement "justification is helpful in believing that 0!=1" for each justification. The aggregate counts for the categories "strongly agree" and "somewhat agree" are 53/62 (85.48%), 49/62 (79.03%), and 44/62 (70.97%) for the three justifications, respectively, and indicate that Justification 1 is the most helpful in believing that zero factorial equals one.

We have to be mindful, however, of the lurking framing effect, which could have made Justification 1 preferred due to its first place in the survey, although we believe that its pronounced simplicity must have made it more attractive to the majority of students. We note that out of 64 students there were at least 62 responses for all statements, reflecting a response rate of 96.88% approximately.

Table 4 shows the counts and percentages for the pre- and post-presentation statements "zero factorial equals one (i.e., 0!=1) is believable", and "overall, the statement that zero factorial equals one (i.e., 0!=1) is more believable after having read the justifications", respectively. The aggregate counts for the categories "strongly agree" and "somewhat agree" are 55/64 (85.94%) for statement 0, and 59/62 (95.15%) for statement 4, an increase of 9.21%. Most notably, the counts for the category "strongly agree" for the pre- and post-presentation statements has an increase of 16.02%. A striking difference observed (in the category "neither agree nor disagree") between the counts for the pre- and post-presentation statements are from 4 to 0 results which in percentages from 6.25% to 0%, respectively. This attests to the scenario of those 4 students who could not conclusively make a decision of their understanding from any of Justifications 1 through 3, are now being able to decide on their preferences from the other four progressively decisive categories. The percentage of students on the fence (neither agree nor disagree), somewhat disagree, or strongly disagree for the pre-presentation statement is 9/64 (14.1%) and 3/62 (4.8%) for the post-presentation statement.

The results of the matched pairs $t$-test was statistically significant ($p$-value = 0.008). The results provide strong evidence that the justifications are helpful in believing that zero factorial equals one.

**Table 1:** Survey response summaries for the statement "*the explanation of symbolic t! is helpful*" (Statements 1a, 2a, and 3a) for the three justifications.

| Response | Justification 1 | | Justification 2 | | Justification 3 | |
|---|---|---|---|---|---|---|
| | Percentage | Count | Percentage | Count | Percentage | Count |
| Strongly agree | 33.87% | 21 | 33.87% | 21 | 32.26% | 20 |
| Somewhat agree | 54.84% | 34 | 48.39% | 30 | 40.32% | 25 |
| Neither agree nor disagree | 8.06% | 5 | 14.52% | 9 | 12.90% | 8 |
| Somewhat disagree | 1.61% | 1 | 1.61% | 1 | 11.29% | 7 |
| Strongly disagree | 1.61% | 1 | 1.61% | 1 | 3.23% | 2 |

**Table 2:** Survey response summaries for the statement "*the pictorial summary by the graph aided further to conceptualize 0! to be 1*" (Statements 1b, 2b and 3b) for the three justifications.

| Response | Justification 1 | | Justification 2 | | Justification 3 | |
|---|---|---|---|---|---|---|
| | Percentage | Count | Percentage | Count | Percentage | Count |
| Strongly agree | 38.71% | 24 | 45.16% | 28 | 43.55% | 27 |
| Somewhat agree | 53.23% | 33 | 46.77% | 29 | 38.71% | 24 |
| Neither agree nor disagree | 4.84% | 3 | 6.45% | 4 | 11.29% | 7 |
| Somewhat disagree | 0.00% | 0 | 0.00% | 0 | 3.23% | 2 |
| Strongly disagree | 3.23% | 2 | 1.61% | 1 | 3.23% | 2 |

**Table 3:** Survey response summaries for the statement "*justification is helpful in believing that* $0! = 1$" (Statements 1c, 2c, and 3c) for the three justifications.

| Response | Justification 1 | | Justification 2 | | Justification 3 | |
|---|---|---|---|---|---|---|
| | Percentage | Count | Percentage | Count | Percentage | Count |
| Strongly agree | 45.16% | 28 | 43.55% | 27 | 35.48% | 22 |
| Somewhat agree | 40.32% | 25 | 35.48% | 22 | 35.48% | 22 |
| Neither agree nor disagree | 8.06% | 5 | 14.52% | 9 | 17.74% | 11 |
| Somewhat disagree | 4.84% | 3 | 4.84% | 3 | 8.06% | 5 |
| Strongly disagree | 1.61% | 1 | 1.61% | 1 | 3.23% | 2 |

**Table 4:** Survey response summaries for the pre- and post-presentation statements "*zero factorial equals one (i.e.,* $0! = 1$*) is believable*" (Statement 0), and "*overall, the statement that zero factorial equals one (i.e.,* $0! = 1$*) is more believable after having read the justifications*" (Statement 4).

| Pre-presentation statement (Statement 0) | | |
|---|---|---|
| Response | Percentage | Count |
| Strongly agree | 46.88% | 30 |
| Somewhat agree | 39.06% | 25 |
| Neither agree nor disagree | 6.25% | 4 |
| Somewhat disagree | 3.13% | 2 |
| Strongly disagree | 4.69% | 3 |

| Post-presentation statement (Statement 4) | | |
|---|---|---|
| Response | Percentage | Count |
| Strongly agree | 62.90% | 39 |
| Somewhat agree | 32.26% | 20 |
| Neither agree nor disagree | 0.00% | 0 |
| Somewhat disagree | 3.23% | 2 |
| Strongly disagree | 1.61% | 1 |

## Discussion

In this paper we have provided three simple and illuminating instructional justifications for zero factorial to be equal to one. In view of the evidence from the survey results, the explanations exhibited the twin goals of addressing an accessible way for all students by revealing that zero factorial is convincingly equal to one while for a range of students providing the rationale for rectifying the perception that zero factorial should be zero. The analysis of the pre- and post-presentation statements have clearly indicated that the justifications are useful in believing that zero factorial equals one.

Overall, the results from the online survey have shown that the students preferred Justification 1. We note that the order of the justifications in which they are presented can be very influential. For instance, the first justification may have convinced many students mathematically and visually that zero factorial equals one, and thus, subsequent justifications scored lower. Another possible reason that Justification 1 was preferred is that it is the most straightforward and may have been the most convincing for students with a weaker mathematical background.

Justification 3 is the best in terms of mathematical rigor, yet it scored lowest in the survey, whereas Justification 1, which scored highest, is the easiest to assimilate but not necessarily the best in terms of mathematical merit. This observation highlights the already made note that the use of the gamma function when teaching statistics at such an early stage is not feasible on pedagogical grounds. We emphasize that teaching introductory level students why zero factorial equals one using the gamma function inhibits the intuitive learning process leading to its conceptual understanding.

In addition to considering the order in which the justifications were presented to the students, our study acknowledges another potential limitation related to instructor bias. As the instructor presenting the justification survey is also a co-author of this paper, there is a possibility of unconscious bias (Kahneman, 2011) influencing the results. To address this concern in future research, employing a different presenter, such as a teaching assistant, to deliver the presentation and administer the online survey could mitigate potential biases (Cole, 1981; StatsCan, 2022). Furthermore, randomly dividing the students into groups and varying the order of the justifications presented to each group would allow for a more robust examination of the impact of presentation sequence (order effects) on student responses.

The appeal of any preferred justification depends on the ability to engage productively based on cognitive skills, acumen and understanding of a given group of students. Having several justifications thus provides a way to meet the educational needs of all students. Inclusive approaches to education have become more prevalent and “psychologists, counsellors, educators, and researchers must critically examine interventions for students with unique educational needs” (Dare et al., 2021). Inclusive education advocates for differentiated instruction, which involves tailoring teaching strategies and materials to meet the individual needs and interests of students (see eg. (Konstantinou-Katzi et al., 2013; Lambert et al., 2021; Roos, 2019)). Educators can provide multiple entry points for understanding zero factorial, offering varying levels of complexity (eg. using the justifications provided in this paper or the gamma function or a combination of the justifications and gamma function using visualizations) and support based on students' prior knowledge, skills, and abilities. This approach allows all students to access the concept at their own pace and level of readiness.

The visual representations accompanying the justifications outlined in this paper can aid in understanding the rationale behind why zero factorial equals one. Visualizations play a vital role in mathematics teaching by providing students with tangible representations of abstract concepts, thereby facilitating deeper understanding and internalization of mathematical ideas (Macnab et al., 2012; Presmeg, 2020). Visualizations engage multiple senses, allowing students to perceive mathematical relationships visually and spatially. By seeing mathematical concepts represented graphically or geometrically, students can make connections between abstract symbols and real-world phenomena, making the concepts more concrete and accessible. Moreover, visual representations often transcend language and cultural barriers, making mathematics more inclusive and accessible to diverse learners.

From the perspective of enriching learning with instructional value driven by survey data, particularly on a special topic, which has a natural counterintuitive connotation, the research undertaken via the three alternative justifications in this article is new to the best of our knowledge.

## CONCLUSION

The justifications presented in this paper may provide instructors alternative ways to initiate students' intuition and hence internalization when teaching the counterintuitive fact that zero factorial equals one. For a diverse group of learners, the justifications have the potential to raise the current level of cognitive skills to inspire differentiated trails of learning (Staples et al., 2012). They offer simple alternative methods, which are mathematically substantiated, necessitating a recognition to an improved conceptual learning. Additionally, it is backed by heuristic evidence in the form of survey results, noting convincingly as to which method is likely to be preferred by introductory level students. In connection to the aforesaid, we have reiterated that students' learning may be facilitated by instructors who apply judiciously at least one of the three justifications in their teaching. Finally, we have illustrated the role of statistical thinking and techniques when making instructional choices.

Moving forward, it is our hope that we continue to explore and develop novel ways to initiate students' intuition, and hence internalization of mathematical concepts like zero factorial. By embracing diverse instructional methods such as incorporating visualizations, we can foster deeper understanding and engagement in mathematics, ultimately enhancing numeracy skills across diverse student populations.

**Acknowledgments:** This version has grown out of numerous discussions with our students in varied classes over several years. We thank the students for their generous feedback on the topic of this paper.

**Ethics:** Ethics approval for this study was obtained from the King's University College Research Ethics Review Committee.

**Author Contribution:** All authors contributed equally to this work.

## REFERENCES

Adams, R., & Essex, C. (2010). *Calculus: A Complete Course* (Ninth ed.). Pearson Canada.

Aunio, P., Korhonen, J., Ragpot, L., Tormanen, M., & Henning, E. (2021). An early numeracy intervention for first-graders at risk for mathematical learning difficulties. *55*, 225-262. https://doi.org/10.1016/j.ecresq.2020.12.002

Bellhouse, D. (1984). A review of optimal designs in survey sampling. *Canadian Journal of Statistics*, *12*(1), 53-65. https://doi.org/10.2307/3314724

Brumwell, S., & MacFarlane, A. (2020). Improving Numeracy Skills of Postsecondary Students: What is the Way Forward? *Higher Education Quality Council of Ontario*. https://heqco.ca/pub/improving-numeracy-skills-of-postsecondary-students-what-is-the-way-forward/

Buglear, J., & Castell, A. (2019). *Stats Means Business: Statistics and Business Analytics for Business, Hospitality and Tourism*. Routledge.

Chen, L., & Zitikis, R. (2017). Measuring and comparing student performance: a new technique for assessing directional associations. *Education Sciences*, *7*(4), 77. https://doi.org/10.3390/educsci7040077

Chen, L., & Zitikis, R. (2020). Quantifying and analyzing nonlinear relationships with a fresh look at a classical dataset of student scores. *Quality & Quantity: International Journal of Methodology*, *54*, 1145-1169. https:// doi.org /10.1007/s11135-020-00979-7

Cole, N. (1981). Bias in Testing. *American Psychologist*, *36*(10), 1067-1077. https://doi.org/10.1037/0003-066X.36.10.1067
Dare, L., Nowicki, E., & Murray, L. (2021). How students conceptualize grade-based acceleration in inclusive settings. *Psychology in the Schools*, *58*(1), 33-50. https://doi.org/10.1002/pits.22435
Devore, J. (2012). *Probability and Statistics for Engineering and the Sciences*. Brooks/Cole, Cengage Learning.
Galligan, L. (2013). A systematic approach to embedding academic numeracy at university. *Higher Education Research & Development*, *32*(5), 734-747. https://doi.org/10.1080/07294360.2013.777037
Gustiningsi, T., Putri, R., Zulkardi, & Hapizah. (2023). Developing a PISA-Like Mathematical Problem: Using Traditional Food Context. *International Journal of Education & Curriculum Application*, *6*(3). https://doi.org/10.31764/ijeca.v6i3.20200
Hairer, E., & Wanner, G. (1996). *Analysis by Its History. Readings in Mathematics*. Springer.
Kahneman, D. (2011). *Thing, Fast and Slow*. Farrar, Straus and Giroux.
Konstantinou-Katzi, P., Tsolaki, E., Meletiou-Mavrotheris, M., & Koutselini, M. (2013). Differentiation of teaching and learning mathematics: an action research study in tertiary education. *International Journal of Mathematical Education in Science and Technology*, *44*(3), 332-349. https://doi.org/10.1080/0020739X.2012.714491
Lambert, R., Imm, K., Schuck, R., Choi, S., & McNiff, A. (2021). “UDL is the What, Design Thinking is the How:” Designing for Differentiation in Mathematics. *23*(3). https://mted.merga.net.au/index.php/mted/article/view/666
Macnab, J., Phillips, L., & Norris, S. (2012). Visualizations and Visualizations in Mathematics Education. *Reading for Evidence and Interpreting Visualizations in Mathematics and Science Education*, 103-145. https://doi.org/10.1007/978-94-6091-924-4_6
Mahmood, M., & Mahmood, I. (2015). A simple demonstration of zero factorial equals one. *International Journal of Mathematical Education in Science and Technology*, *47*(6), 959-960. https://doi.org/10.1080/0020739X.2015.1119896
Mahmood, M., & Romdhane, M. (2015). An algebraic justification that zero factorial equals one. *The Mathematical Scientist*, *40*, 134-135. https://www.appliedprobability.org/publications/the-mathematical-scientist
Meaders, C., Lane, A., Morozov, A., Shuman, J., Toth, E., Stains, M., & al., e. (2020). Undergraduate student concerns in introductory STEM courses: what they are, how they change, and what influences them. *Journal for STEM Education Research*, *2*, 195-216. https://doi.org/10.1007/s41979-020-00031-1
Murray, L., & Wilson, J. (2021). Generating data sets for teaching the importance of regression analysis. *Decision Sciences: Journal of Innovative Education*, *19*(2), 157-166. https://doi.org/10.1111/dsji.12233
NACOME. (2023). The 2nd National Conference on Mathematics Education (NACOME) 2021: Mathematical Proof as a Tool for Learning Mathematics. *2811*. https://doi.org/10.1063/12.0016981
Presmeg, N. (2020). Visualization and Learning in Mathematics Education. *Encyclopedia of Mathematics Education*. https://doi.org/10.1007/978-3-030-15789-0_161
Putra, A., Zulkardi, Putri, R., Susanti, E., & Nusantara, D. (2023). Teaching the Rule of Product using Nature Tourism Routes. *Tarbawi: Jurnal Iimu Pendidkan*, *19*(2). https://doi.org/10.32939/tarbawi.v19i2.3168
Rawani, D., Putri, R., Zulkardi, & Susanti, E. (2023). Explore the South Sumatera's traditional dance as the concept to support mathematical learning. *AIP Conference Proceedings*, *2811*(1). https://doi.org/10.1063/5.0142305

Risdiyanti, I., Zulkardi, Putri, R., Prahmana, R., & Nusantara, D. (2024). Ratio and proportion through realistic mathematics education and pendidikan matematika realistik Indonesia approach: A systematic literature review. *Jurnal Elemen*, *10*(1). https://doi.org/10.29408/jel.v10i1.24445

Roberts, T., Jackson, C., Mohr-Schroeder, M., Bush, S., Maiorca, C., Cavalcanti, M., & al., e. (2018). Students' perceptions of STEM learning after participating in a summer informal learning experience. *International Journal of STEM Education*, *5*. https://stemeducationjournal.springeropen.com/articles/10.1186/s40594-018-0133-4

Roos, H. (2019). Inclusion in mathematics education: an ideology, a way of teaching, or both? *Education Studies in Mathematics*, *100*, 25-41. https://doi.org/10.1007/s10649-018-9854-z

Staples, M., Bartlo, J., & Thanheiser, E. (2012). Justification as a teaching and learning practice: Its (potential) multifaceted role in middle grades mathematics classrooms. *The Journal of Mathematical Behavior*, *31*. https://doi.org/10.1016/j.jmathb.2012.07.001

StatsCan. (2022). Statistics 101: Statistical Bias. *Statistics Canada*. https://www.statcan.gc.ca/en/wtc/data-literacy/catalogue/892000062022005

Thompson, M. (1997). *Theory of Sample Surveys*. Chapman and Hall.

Wiggins, B., Eddy, S., Wener-Fligner, L., Freisem, K., Grunspan, D., Theobald, E., & al., e. (2017). ASPECT: A survey to assess student perspective of engagement in an active-learning classroom. *Life Sciences Education*, *6*. https://doi.org/10.1187/cbe.16-08-0244

Wu, C., & Thompson, M. (2020). *Sampling Theory and Practice*. Springer.

Zhai, X., Gu, J., Liu, H., Liang, J., & Tsai, C. (2017). An experiential learning perspective on students' satisfaction model in a flipped classroom context. *Journal of Educational Technology & Society*, 198-210. http://www.jstor.org/stable/jeductechsoci.20.1.198